\documentclass{amsart}
\usepackage{lmodern}
\usepackage{textcomp}
\usepackage{amsfonts,amssymb,amsmath,amscd,amstext}
\usepackage{esint}
\usepackage{enumerate}
\usepackage{cancel,cite}
\usepackage{bbm}
\usepackage[pdftex]{graphicx}
\usepackage{floatflt}
\usepackage{mathrsfs}
\usepackage{ulem,lipsum}
\usepackage{mathdots}

\usepackage[dvipsnames]{xcolor}
\usepackage{hyperref}
\hypersetup{
	colorlinks=true,
	linkcolor=blue,
	urlcolor=black, 
	citecolor=blue
}
\usepackage[noabbrev,capitalize,nameinlink]{cleveref}
\usepackage[margin=1.45in]{geometry} 

\allowdisplaybreaks

\def\vs{\vspace{1mm}}

\def\dx{\,{\rm d}x}

\def\dv{\,{\rm d}v}
\def\dw{\,{\rm d}w}
\def\dt{\,{\rm d}t}

\def \d{\,{\rm d}}

\def\er{\mathbb R}

\def \K {{K}}
\def \W {\mathcal{W}}
\def \Wc {\mathcal{W}}
\def \Ec {\mathcal{E}}

\def \R {\mathbb{R}}

\def \and {$\, \times\,$}

\makeatother

\def\mean#1{\mathchoice%
	{\mathop{\kern 0.2em\vrule width 0.6em height 0.69678ex depth -0.58065ex
			\kern -0.8em \intop}\nolimits_{\kern -0.4em#1}}%
	{\mathop{\kern 0.1em\vrule width 0.5em height 0.69678ex depth -0.60387ex
			\kern -0.6em \intop}\nolimits_{#1}}%
	{\mathop{\kern 0.1em\vrule width 0.5em height 0.69678ex
			depth -0.60387ex
			\kern -0.6em \intop}\nolimits_{#1}}%
	{\mathop{\kern 0.1em\vrule width 0.5em height 0.69678ex depth -0.60387ex
			\kern -0.6em \intop}\nolimits_{#1}}}
\newcommand{\vertiii}[1]{{\left\vert\kern-0.25ex\left\vert\kern-0.25ex\left\vert #1 
		\right\vert\kern-0.25ex\right\vert\kern-0.25ex\right\vert}}
\makeatletter        
\newcommand{\linethrough}{\mathpalette\@thickbar}
\newcommand{\@thickbar}[2]{{#1\mkern0mu\vbox{
    \sbox\z@{$#1#2\mkern-0.5mu$}%
    \dimen@=\dimexpr\ht\tw@-\ht\z@+2\p@\relax % The +2 represents the vertical shift of the line.
    \hrule\@height0.5\p@ % The 0.5 represent the thickness on the line.
    \vskip\dimen@
    \box\z@}}}
\makeatother
\newcommand{\mathstrike}[1]{\ensuremath{\linethrough{#1}}}
\newcommand{\nra}[1]{\mathstrike{\lVert} #1 \rVert}

\allowdisplaybreaks
\newtheorem{theorem}{Theorem}[section]

\newtheorem{defn}[theorem]{Definition}
\newtheorem{rem}[theorem]{Remark}
\newtheorem{lemma}[theorem]{Lemma}
\newtheorem{ass}[theorem]{Assumption}

\numberwithin{equation}{section}

\newcommand{\ern}{\mathbb{R}^n}
\newcommand\eps\varepsilon

\newcommand{\rr}{\varrho}

\newcommand{\norma}[1]{{\left\|#1\right\|}}
\newcommand{\snr}[1]{\lvert #1\rvert}
\newcommand{\Lc}{\mathcal{L}}
\newcommand{\tail}{\textup{Tail}}

\newcommand{\F}{\varphi}

\renewcommand{\Sigma}{\varSigma}
\renewcommand{\Lambda}{\varLambda}
\renewcommand{\Psi}{\varPsi}
\renewcommand{\Phi}{\varPhi}

\title[Boundedness for nonlinear nonlocal Fokker-Planck equations]{Boundedness estimates		
for nonlinear nonlocal kinetic
Kolmogorov-Fokker-Planck equations}

\author[F. Anceschi]{Francesca Anceschi}  
\address[F. Anceschi]{Dipartimento di Ingegneria Industriale e Scienze Matematiche
	\newline\indent Universit\`a Politecnica delle Marche
	\newline\indent Via Brecce Bianche, 12, 60131 Ancona, Italy} \email{\url{f.anceschi@staff.univpm.it}}

\author[M. Piccinini]{Mirco Piccinini} 
\address[M. Piccinini]{Dipartimento di Matematica
	\newline\indent Universit\`a di Pisa
	\newline\indent L.go~B.~Pontecorvo~5, 56127, Pisa, Italy}
\email{\url{mirco.piccinini@dm.unipi.it}}

\newlength{\defbaselineskip}
\begin{document}

	\keywords{Kolmogorov-Fokker-Planck equations, kinetic equations, fractional Sobolev spaces, fractional Laplacian, nonlinear operators\vspace{1mm}}
	
	\subjclass{35Q84, 35B45,  35B65, 47G20,  35R11, 35R05\vspace{1mm}}

	\begin{abstract}
		We investigate local regularity properties of weak solutions to a broad class of nonlinear 
        nonlocal kinetic Kolmogorov-Fokker-Planck equations.        
    In particular, we focus on proving an 
		interpolative apriori boundedness estimate  for weak 
        		subsolutions in terms of a tail term encoding the nonlocal contributions of the diffusion.
	\end{abstract}
    
%%% Title
\maketitle

%----------------------------------------------------------------------
\section{Introduction}
In this work we deal with a wide class of kinetic equations, whose diffusion term is driven by an integro-differential operator of differentiability 
order~$s \in (0,1)$, which is allowed to be nonlinear with at most  quadratic growth. More specifically, we investigate local properties of 
weak solutions~$f\equiv f(t,x,v)$ to the following class of equations
	\begin{equation}\label{pbm}
		(\partial_t +v \cdot \nabla_x) f (t,x,v) =
		  \Lc f(t,x,v)
		  \qquad \text{for}~(t,x,v) \in \er \times \ern \times \ern\,,
	\end{equation}
	where the nonlocal %{term} 
    operator $\Lc$ is given by
\begin{equation}\label{eq:main-op}
	\Lc f(t,x,v) := P.~\!V. \int_{\ern} \varPhi\left(\frac{f(t,x,w)-f(t,x,v)}{\snr{v-w}^s}\right)\frac{\dw}{\snr{v-w}^{n+s}}.
\end{equation}
Here, the symbol $P.~\!V.$ stands for “in the principal value sense”,~$s \in (0,1)$ and the nonlinearity~$\varPhi$ satisfies the following assumption.

\begin{ass}\label{assumption}
We assume that~$\varPhi : \er \to \er$ is an odd function such that for any~$\tau,\tau' \in \er$ and some~$\Lambda \geq 1$, it holds
\begin{align}\label{eq:nonlinearity-ellipt-1}
    &\snr{\varPhi(\tau)-\varPhi(\tau')} \leq \Lambda \snr{\tau-\tau'} \quad \\ \nonumber
    & \qquad \text{and} \quad \big(\varPhi(\tau)-\varPhi(\tau')\big)(\tau-\tau') \geq \Lambda^{-1}\snr{\tau-\tau'}^2.
\end{align}
\end{ass}
\noindent
Note that, since $\Phi$ is an odd function, it is also true that $\Phi (0)=0$.

 \vspace{2mm}
 As a prototype for Equation~\eqref{pbm}, even though in this scenario the difficulties arising when dealing with a nonlinear operator vanish, one can consider~$\Phi$ to be the identity. Then, the diffusion
 in velocity coincides with the classical 
 fractional Laplacian
 and in this setting Equation~\eqref{pbm} 
 does reduce to
\begin{equation}\label{eq:model}
	(\partial_t +v \cdot \nabla_x) f +  (-\Delta_v)^sf=0,
\end{equation}
Equation \eqref{pbm} can be seen as a nonlocal extension of the nonlinear (local) equation studied in \cite{Nystrom} by Nystr\"om and Garain where the vector filed 
$A(\cdot)$, generating the nonlinear structure of the equation, satisfies standard ellipticity and quadratic growth assumptions, %and 
which can be recovered, say by taking $s=1$ in \eqref{pbm}. 
Aside from the aforementioned paper, where interior 
regularity {\it \'a la} De Giorgi-Nash-Moser is addressed we recall
\cite{LascialfariMorbidelli} by Lascialfari and Morbidelli, where the 
well-posedness of a Dirichlet problem in the local quasilinear case is addressed, as well as the recent \cite{KLN25} where precise pointwise estimates in the spirit of nonlinear potential theory and fine gradient regularity results under borderline assumptions on the data were achieved in a setting analogous to \cite{Nystrom}.

\vspace{1mm}

For what regards the linear nonlocal case, the study of the regularity  theory and qualitative properties of fractional kinetic
 equations has recently witnessed a substantial growth by attracting the attention of different mathematical communities.
This is partially due to the appearence of nonlocal kinetic equations in several, even seemingly unrelated, models,  as, e.\ \!g., in Finance, in order to describe the evolution of 
Asian options, where the drift term is connected with risk-free interest rates as well as in Gas Dynamics where they appear as linearized models for the Boltzmann equation 
without cutoff. In this scenario, a priori boundedness and further regularity estimates would be very useful results in order to tackle with well-posedness issues 
and long-time behavior studies; see for instance the famous result on the trend to global equilibrium of Desvillettes and Villani \cite{Villani} for the Boltzmann 
equation without cutoff where solutions are assumed apriori to be smooth up to the boundary. In this scenario,
the weak regularity theory for nonlocal equations 
has been the main focus of various recent  efforts by different communities; see~\cite{APR23} and the references therein. In particular, we refer the reader to the H\"older 
regularity results in~\cite{Sto19}, possibly including unbounded source terms, as well as the ones in~\cite{Loh22} covering more general, possibly nonsymmetric diffusion 
operators. Furthermore, regarding classical estimates, we mention the very recent breakthrough counterexample to the  classical Harnack inequality~\cite{KW-controesempio}, as well as its related new formulation in~\cite{APP25}, where a strong Harnack inequality with tail is proved provided that solutions have $\sigma$-summable 
nonlocal tail along the transport variables for some~$\sigma > \sigma^\star(n,s)$, which is in fact naturally implied by the usual assumptions considered in literature, e.~\!g., from the 
usual mass density boundedness (as for the Boltzmann equation without cut-off), and in clear accordance with the aforementioned counterexample in~\cite{KW-controesempio}. Still in the flavor of Harnack-type inequalities, it is worth mentioning the very recent paper~\cite{Loh24c}, in which amongst other results, the author proves a 
strong Harnack inequality for {\it global solutions}, a priori bounded, periodic in the space variable, and under an integral monotonicity-in-time assumption (see~Definition~2.2 
there). Finally, we mention~\cite{WD17} for the proof of the existence of weak solutions, and~\cite{Imb05} for existence, uniqueness and regularity results for solutions in the 
viscosity sense to fractional linear kinetic equations. Always regarding these existence and uniqueness issues, we also recall the very recent works~\cite{niebel1, niebel2}. 

\vspace{2mm}

For what regards more general nonlinear nonlocal kinetic equations
to the best of our knowledge, our contribution would be the first.
In this respect, forthcoming Theorem~\ref{loc bdd} serves as a first step in the direction of proving that solutions to~\eqref{pbm} enjoy classical qualitative properties and 
extend the results already available in the local case proved in the aforementioned \cite{Nystrom}; 
see Section \ref{sec:nonlinear} for further information on other types of growth (subquadratic, or superquadratic ones)
in the same flavour of~\cite{APT22,Liao22,Tav24}. 
    
    Aside from the novelty of our results, these quantitative estimates are very useful when dealing with local regularity, or qualitative properties of solutions to~\eqref{pbm}. However, even proving a $L^2$-$L^\infty$ estimate for nonlinear nonlocal kinetic equations is not a simple task. Indeed, even in the %nonlocal
    linear case -- as proven in the aforementioned work~\cite{KW-controesempio} -- it is not possible, in general, to bound the $L^\infty$ norm of a solution in terms of only local quantities {\it even starting from globally bounded solutions}. Moreover, a deeper analysis of the counterexample in~\cite{KW-controesempio} shows that  such supremum estimate remains false also when an error term is added on its right-hand side -- basically a tail-type contribution as in~\eqref{eq:tail} -- if the tail belongs to~$L^\sigma$, for~$\sigma < (n(1+2s))/(2s)$.
    
In order to balance in a quantitative way the nonlocal behavior of the diffusion in velocity with the lack of ellipticity in the spacial variable given by the additional transport term, we have to work under a sufficient integrability assumption on the {\it the nonlocal tail of a function}~(\cite{DKP14,DKP16}) defined as
    \begin{equation}\label{eq:tail}
    \tail(f;B_r(v_{\rm o})):= r^{2s}\int_{\ern \setminus B_r(v_{\rm o})} \frac{\snr{f(t,x,v)}}{\snr{v_{\rm o}-v}^{n+2s}}\dv\,.
    \end{equation}	
Indeed, the finiteness of the $L^\sigma$-energy of the tail term is a turning point in the local analysis of~\eqref{pbm}, and appears to be in contrast with most of the parabolic 
literature, where nonlocal effects have been compensated via a (sharp)  $L^1$-tail (see \cite{KW23}), which however is critical with respect to kinetic scalings.
Moreover, even if by definition weak solutions are not required to have finite $L^\sigma$-tail, the usual constraints on the mass
observable, see \cite{IS20}, plainly imply our requirements on the $L^\sigma$-energy of the nonlocal tail.
    
Lastly, as in~\cite{APP25}, the backbone of the proof of a $L^2$-$L^\infty$ estimate is a hypoelliptic gain of integrability, which is proven by making use of the fundamental solution of the fractional Kolmogorov equation. More specifically, as in the classical framework for kinetic equations~(\cite{PP04}), the transfer of regularity is based on treating as source term the difference between the constant coefficients diffusion operator and the one with measurable entries, and then estimating its $L^2$-norm tracking down the long-range interactions appearing as $L^\sigma$-norm of the tail quantity~\eqref{eq:tail} on the right-hand side; see, in particular,~\cite[Theorem 1.4]{APP25}. However, as well as for velocity averaging lemmas, such a procedure can not be pursued in a very general nonlinear setting. Hence, we focus on a fractional nonlinear case with quadratic growth, and we prove that subsolutions to \eqref{pbm} satisfy interpolative estimates in terms of their local and nonlocal contributions.
    The interpolative nature of the estimates below lies specifically in the arbitrariness in the choice of the parameter~$\delta$, which plays the role of an interpolation coefficient between the local and nonlocal part of the estimate. We also remark that essentially the~$\tail(\cdot)$ behaves as a source term. Hence, in this respect the forthcoming lower bound on the integrability condition in~\eqref{e-star} is the expected one also with respect to the local analogue of Kolmogorov equation; see \cite{AR22,Hou25}.

%------------------------------------------------------------------------
\begin{theorem}\label{loc bdd}
      	\label{thm_bdd}
		Let~$\Omega:= (t_1,t_2)\times \Omega_x \times \Omega_v \subset \er^{1+2n}$ be a domain,~$s \in (0,1)$ and let~$N_{s}$ be the homogeneous dimension in~\eqref{def:homo-dim}.  Assume that~$f \in \Wc$ is a weak  subsolution to
	\eqref{pbm}.
	If~$\tail(f_+;B) \in L^\sigma_{\rm{loc}}((t_1,t_2)\times\Omega_{x})$ for any~$B \Subset \Omega_v$, for some~$\sigma$ satisfying~
	\begin{equation}\label{e-star}
		\sigma > \frac{N_{s}}{2s}\,,
	\end{equation}
	then, for any~$ {Q}_{r}(z_{\rm o})  \subset Q_1(z_{\rm o}) \subset  \Omega$ and any~$\delta \in (0,1]$,  it holds
	\begin{eqnarray*}
		\sup_{{Q}_{\frac{r}{2}}(z_{\rm o})}f & \, \leq \, & c\,\left(\frac{\langle v_{\rm o}\rangle^{\frac{3}{2}}}{ r^{\frac{n}{\sigma}+2(n+2s)}\delta }\right)^\frac{\sigma N_{s}}{2s\sigma-N_{s}}
		\left( \mean{Q_{r}(z_{\rm o})}  \snr{f}^2 \, \dt \dx \dv  \right)^{\frac12}  \\
		&&+
		 \ \delta  \left( \mean{U_{r}(t_{\rm o},x_{\rm o})} \snr{ \tail(f_+;B_\frac{r}{2}(v_{\rm o}))}^{\sigma} \dt \, \dx \right)^{\frac1p}  \,, \notag 
	\end{eqnarray*}
	where~$c \equiv c(n,s,\Lambda,\sigma)>0$.
\end{theorem}
\begin{rem}{\rm 
    It is possible to extend this results to an even wider class of nonlinear nonlocal kinetic operators defined as 
    \begin{align} \label{pbm1}
        \partial_tf(t,x,v) + \nabla_x \cdot (v f(t,x,v) + u(t,x) f(t,x,v)) =
		  \Lc f(t,x,v) 
    \end{align}
    for $(t,x,v) \in \er \times \ern \times \ern$, 
    where $\Lc$ is defined in \eqref{eq:main-op} and $u: \mathbb{R} \times \mathbb{R}^n \to \mathbb{R}^n$ is a vector field. 
    This class is the nonlinear nonlocal counterpart of the one considered in
    \cite{ZhuAv}, and to extend our analysis to this case one needs to follow the theory presented in
    the aforementioned paper, since the major difference between \eqref{pbm} and \eqref{pbm1} lies in the structure of the transport, where a vector field not depending on $v$ is added. On the other hand, a more careful investigation would be needed 
    in order to consider a nonlinear nonlocal equation of the following type 
    \begin{align*}% \label{pbm1}
        (\partial_t + b(v) \cdot \nabla_x) f(t,x,v) =
		  \Lc f(t,x,v) 
    \end{align*}
    for $(t,x,v) \in \er \times \ern \times \ern$, 
    where $\Lc$ is defined in \eqref{eq:main-op} and $b$ satisfies a non-degeneracy assumptions, see for instance Equation (1.3) of 
    \cite{Zhub}, where the analogous analysis in the linear local case is carried out.}
\end{rem}

%------------------------------------------------------------------------
The proof of the above theorem strongly relies on the combination of a Caccioppoli inequality, see forthcoming Lemma \ref{thm ccp}, together with a higher integrability result for subsolutions to \eqref{pbm} proved by making use of the fundamental solution of the linearized version of equation \eqref{pbm} and following the approach presented in \cite{APP25}.

\vspace{2mm}
\mbox{}
\\ {\bf Outline  of the paper.} In Section~\ref{sec:preliminaries} we introduce preliminary notions about the functional and geometrical setting 
   of this work. Section~\ref{sec_cc} is devoted to the proof of  Theorem~\ref{thm_bdd}. 
   Section~\ref{sec:nonlinear} we discuss in detail the challenges one needs to face when dealing with the  nonlocal $p$-growth case.
 
 \vspace{2mm}
 \mbox{}
 \\ {\bf Aknowledgements.} FA is partially supported by the INdAM - GNAMPA project ``Variational problems for Kolmogorov equations: long-time analysis and regularity estimates'',  CUP\_E55F22000270001. MP is {partially} supported by the INdAM - GNAMPA project 
 	``Fenomeni non locali in problemi locali", CUP\_E55F22000270001.
 	Both authors are partially supported by the INdAM-GNAMPA Project ``Problemi non locali: teoria cinetica e non uniforme ellitticità'', CUP\_E53C22001930001.
 
\section{Notation and preliminaries} \label{sec:preliminaries}
In this Section, we recall some known results about our underlying geometrical and functional setting. 
After fixing the notation, we introduce an appropriate geometric framework to study integral kinetic equations; then, we recall some properties of fractional Sobolev spaces and the functional setting required to deal with Equation~\eqref{pbm}.

\subsection{Notation and geometric framework}
We denote with~$c$ a positive universal constant greater than one, which may change from line to line. For the sake of readability, dependencies of the constants will be often omitted within the chains of estimates, therefore stated after the estimate. Relevant dependencies on parameters will be emphasized by using parentheses. For any~$\mathcal{O} \subset \er^n$ we denote with~$\mathbf{1}_{\mathcal{O}}$ the indicator function of~$\mathcal{O}$.
As customary, for any~$r>0$ and any~$y_{\rm o} \in \er^n$ we denote by~$
B_r(y_{\rm o}) \equiv B(y_{\rm o};r):=\big\{y \in \er^n \,:\, \snr{y-y_{\rm o}}< r\big\}$\,,
the open ball with radius~$r$ and center~$y_{\rm o}$. We shall often abbreviate $B_1 \equiv B_1(0)$, where we denote with $0_{\ern}:= 0$. 
For any measurable function~$g$, we define  the positive 
and negative part of~$g$ as~$g_\pm(y):=\max\{\pm g(y),0\}$.  We denote with $z := (t,x,v) $ a generic point of $\er^{1+2n}$. We shall often abbreviate $B_1 \equiv B_1(0)$, where we denote with $0_{\ern}:= 0$. 
For any measurable function~$g$, we define  the positive 
and negative part of~$g$ as~$g_\pm(y):=\max\{\pm g(y),0\}$. Given any open set~$\mathcal{O} \subset \er^{1+2n}$, with positive Lebsegue measure $\snr{\mathcal{O}}>0$ we denote with
\begin{align} \label{nra}
\nra{g}_{L^p(\mathcal{O})} := \left(\mean{\mathcal{O}}\snr{g}^p \d z \right)^\frac{1}{p}.
\end{align}
In a similar fashion, as for the geometrical setting of the Boltzmann kernel~\cite{IS} or as in~\cite{MPP24}, we start by endowing~$\er^{1+2n}$ with the following Galilean transformation
\begin{equation}\label{eq:group-law}
	z_{\rm o} \circ z : = (t+t_{\rm o},\, x+x_{\rm o}+tv_{\rm o}, \, v+v_{\rm o}).
\end{equation}
With respect to the group law~$\circ$, the couple~$(\er^{1+2n}, \,\circ)$ is a Lie group with identity element~$0:=(0,0,0) \in \R^{1+2n}$ and inverse element~
\begin{equation*}
	z^{-1} = (t,x,v)^{-1}= (-t,\, -x+tv,\, -v) \qquad \text{for any } \, (t,x,v) \in \er^{1+2n}.
\end{equation*} 

Furthermore, for any~$r>0$, we consider the usual fractional nonlinear kinetic scaling~$D(r):  \er^{1+2n} \mapsto \er^{1+2n}$ defined by
\begin{equation}\label{eq:dil}
	\delta_r (t,x,v):=(r^{2s}t,\, r^{1+2s}x,\, r v). 
\end{equation}
Then for any~$r>0$, we denote by~${Q}_r$ a cylinder centered in the origin of radius~$r$; that is,
\begin{align} \label{eq:cyl}
{Q}_r \equiv {Q}_r({0}):= U_r(0,0)\times B_r(0) = (-r^{2s},0]\times B_{r^{1+2s}}(0)\times B_r(0)\,.
\end{align}
For every~$z_{\rm o} \in \er^{1+2n}$ and for every~$r>0$, the {\it slanted} cylinder~${Q}_r(z_{\rm o})$ is defined as follows,
\begin{eqnarray*}
	&& {Q}_r(z_{\rm o})\!\!\! \ :=
	\big \{ z:=(t,x,v) \in \er^{1+2n}: \, -r^{2s} < t -  t_{\rm o} \leq  0 , \\* 
&&\hspace{3,5cm}   \ | x - x_{\rm o} - (t - t_{\rm o}) v_{\rm o}| <r^{1 + 2s},\ \snr{ v- v_{\rm o} } < r \big\}. \notag
\end{eqnarray*}

We denote with~$N_{s}$ the {\it homogeneous dimension} related to~\eqref{eq:dil} defined as
\begin{equation}\label{def:homo-dim}
N_{s}:=  n(2+2s)+2s.
\end{equation}	
Such quantity encodes the scaling properties of the underlying kinetic scalings. Indeed, we have that~$|Q_r| = r^{N_{s}}|Q_1|$, and in general~$|\delta_r(\Omega)|= r^{N_{s}}|\Omega|$, for any Lebesgue measurable sets~$\Omega\subset \er^{1+2n}$.

\subsection{The functional framework}
Lastly, we now introduce the family of related function spaces. For~$s\in (0,1)$  denote with~$H^s(\mathcal{O})$ the classical fractional Sobolev space
\[
H^s(\mathcal{O}) := \left\{f \in L^2(\mathcal{O}) \; : \; \left[ f \right]_{H^s(\mathcal{O})}  < +\infty \right\},
\]
where the fractional seminorm~$\left[ f \right]_{H^s(\mathcal{O})} $ is the usual one  defined via the Gagliardo kernel
\[
\left[ f \right]_{H^s(\mathcal{O})}  := \left(\iint_{\mathcal{O}\times \mathcal{O}} \frac{\snr{f(v)-f(w)}^2}{\snr{v-w}^{n+2s}} \, \dv\dw\right)^{1/2}\,,
\]
and where we have equipped~$H^s$ with the usual norm
\[
\norma{f}_{H^s(\mathcal{O})} := \norma{f}_{L^2(\mathcal{O})} + \left[ f \right]_{H^s(\mathcal{O})}.
\]
A function~$ f$ belongs to $H_{\rm loc}^s(\mathcal{O}))$ if~$f\in H^s(\mathcal{O})')$ whenever~$\mathcal{O})' \Subset \mathcal{O})$. We will denote with~$H^{-s}(\ern)$ the dual of~$H^s(\ern)$ and denote with~$\langle \cdot,\cdot \rangle_{H^{-s},H^s}$ the usual duality paring between~$H^{-s}$ and~$H^s$.
Let us remark that, via Riesz-Fr\'echet's Representation Theorem 
for any~$f \in H^{-s}(\ern)$, there exists two functions~$h_{\rm o}$,~$h_1 \in L^2(\ern)$ such that
$f = h_1 + (-\Delta_v)^{s/2}h_{\rm o}$ {and} $\|h_{\rm o}\|_{L^2(\ern)} + \|h_1\|_{L^2(\ern)} \approx \|f\|_{H^{-s}(\ern)}$.
For any~$f \in H^s(\ern)$ we define~$\Lc f$ as an element of~$H^{-s}(\ern)$ that acts on~$\phi \in H^s(\ern)$ via
\[
\langle \Lc f , \phi \rangle_{H^{-s},H^s} = \frac{1}{2}\iint_{\ern\times\ern}\varPhi\left(\frac{f(v)-f(w)}{\snr{v-w}^s}\right)\frac{(\phi(v)-\phi(w))}{\snr{v-w}^{n+2}}\dw\dv.
\]

Consider the following  tail space
\[
L^1_{2s}(\ern):= \left\{g \in L^1_{\rm{loc}}(\ern)\, : \,  \int_{\ern}\frac{|g(v)|}{(1+\snr{v})^{n+2s}}\dv < \infty \right\},
\]
as firstly defined in~\cite{KKP16}. Then, given~$\Omega:=    (t_1,t_2)\times \Omega_x \times \Omega_v \subset \er^{1+2n}$,  
we denote by~$\Wc$ the natural functions space to which weak solutions to~\eqref{pbm} belong to, and defined as
 \begin{multline*} 
 	\Wc := \Big\{ f \in L^2_{\textrm{loc}}((t_1, t_2)\times \Omega_x ;\, H^s_{\textrm{loc}}(\Omega_v))\cap L^1_{\textrm{loc}}((t_1, t_2)\times \Omega_x ;L^1_{2s}(\ern)) 
 	\\*: (\partial_t +v \cdot \nabla_x)f \in {L^2_{\textrm{loc}} ((t_1,t_2)\times \Omega_x;\,  H^{-s}(\ern )} \Big \}.
 \end{multline*}

We are now in a position to recall the definition of weak sub- and supersolution.
\begin{defn}\label{weak-sol-int}
	A function~$f \in \Wc$ is a \,{\rm weak subsolution} {\rm (}resp.,~{\rm supersolution}{\rm)} to~\eqref{pbm} in~$\Omega$ if 
	\begin{align*}
		&\int_{t_1}^{t_2}\int_{\Omega_x} \iint_{\mathbb{R}^n \times \mathbb{R}^n} \varPhi\left(\frac{f(t,x,v)-f(t,x,w)}{\snr{v-w}^{s}}\right)
		\frac{\psi(t,x,v)-\psi(t,x,w)}{\snr{v-w}^{n+s}}\dw\dv \dx\dt\notag\\
        & \quad - \,\int_{t_1}^{t_2}\int_{\Omega_x}
		\int_{\Omega_v}  \,f (t,x,v) \,  (\partial_t+v	\cdot \nabla_x)  \psi(t,x,v) \dv \dx \dt
		\, \, \leq\, 0 \, \, \big(\ge\,0, \, \text{resp.}\big),
	\end{align*}
	for any nonnegative~$\psi \in L^2  ((t_1,t_2)\times \Omega_x ;H^s_0(\Omega_v))$.
	A function~$f \in \Wc$ is a {\rm weak solution} to~\eqref{pbm} if it is both a weak sub- and supersolution.
\end{defn}

\section{The $L^2$-$L^\infty$ estimate}\label{sec_cc}

From now on, for any two functions $f,g \in H^s(\mathbb{R}^n)$, and $t \in \mathbb{R}$ and $x \in \mathbb{R}^n$, we define the bilinear operator 
\begin{align*}
    \mathcal{E}(f, g) (t,x):= \iint_{\mathbb{R}^n \times \mathbb{R}^n} \varPhi\left(\frac{f(t,x,w)-f(t,x,v)}{\snr{v-w}^{s}}\right)\frac{g(w)-g(v)}{\snr{v-w}^{n+s}}\dw\dv.
\end{align*}
\noindent

As in the classical theory, a fractional Caccioppoli-type inequality is needed in order to built the proper iteration scheme. This extends to the nonlinear setting the approach seen in~\cite[Lemma 3.1]{APP25}; see in particular Step~1 there. 

\begin{lemma}[Caccioppoli inequality]\label{thm ccp}
   	Let~$s \in (0,1)$ and $\Omega\subset \er^{1+2n}$ be a bounded domain.
   Let~$f $ be a weak subsolution to~\eqref{pbm} in~$\Omega$ according to {\rm Definition~\ref{weak-sol-int}}. For any $r \in (0,1)$ such that~$Q_{r}(z_{\rm o}) \subset \Omega$ the following estimate holds true for any $\kappa \in \er$, $\sigma>2$ and $\rr \in (0,r)$
   	\begin{eqnarray} \label{caccioppoli-def}
   	 	&& \sup_{t \in (-\rr^{2s}+t_{\rm o},t_{\rm o}]}\int_{Q_\rr^t(z_{\rm o})} (f-\kappa)^2_+\dx\dv+  \int_{U_\rr(t_{\rm o},x_{\rm o})} [(f-\kappa)_+]^2_{H^s(B_\rr(v_{\rm o}))} \dt \dx \nonumber \\*[0.5ex]
   		&& \quad\leq   \frac{c\langle v_{\rm o}\rangle} 
              {(r-\rr)^{2(1+s)}}   \int_{Q_r(z_{\rm o})} (f-\kappa)_+^2\dt \dx \dv \\
        &&\qquad +\frac{c\,|Q_r(z_{\rm o})\cap \{f > \kappa\}|^{\frac12-\frac{1}{\sigma}}}{({r}-{\varrho})^{2(n+2s)}}
       \|(f-\kappa)_+\|_{L^2(Q_r(z_{\rm o}))} \|\tail((f-\kappa)_+;B_{{r}}(v_{\rm o})) \|_{L^\sigma \left( U_r(t_{\rm o},x_{\rm o})\right) } \,, \nonumber
   	\end{eqnarray}
	where the constant~$c>0$ depends only on $\Lambda$, $n$ and $s$.
\end{lemma}

It is possible to prove an analogous result even for nonquadratic growths. For further information on this matter, we refer the reader to upcoming Section \ref{sec:nonlinear}.

\begin{proof}
 Let us fix $r \in (0,1)$ such that~$Q_{r}(z_{\rm o}) \subset \Omega$. Then, with no loss of generality we assume that the main cylinder is centered at the origin.
Indeed, by~\cite[Lemma~5.1]{Sto19} the function~$\tilde{f}(z):=f(z_{\rm o}\circ z)$ satisfies
\[
\partial_t  \tilde{f}  + v\cdot \nabla_x\tilde{f} = \Lc \tilde{f}  \quad \text{in}~z_{\rm o}^{-1}\circ \Omega.
\]
Fix~$0 < \varrho < r <1$ such that~$Q_{r} \equiv Q_{r}(0) \subset z_{\rm o}^{-1}\circ \Omega$, and consider a weak subsolution to \eqref{pbm}. Then, we 
consider two function that for $-r^{2s}< \tau_1 < \tau_2 <0$ and $0 < \delta < \min \{ \tau_1 + r^{2s}, - \tau_2 \} $
are defined as
     \[
     \theta_\delta(t) :=
     \begin{cases}
         0 & \quad \text{for} \ -r^{2s} \leq t \leq \tau_1 -\delta,\\
         1+\frac{t-\tau_1}{\delta} & \quad \text{for} \ \tau_1 -\delta < t \leq \tau_1,\\
         1 & \quad \text{for} \ \tau_1 < t \leq \tau_2,\\
         1- \frac{t-\tau_2}{\delta} & \quad \text{for} \ \tau_2 < t \leq \tau_2+\delta,\\
         0 & \quad \text{for} \ \tau_2 +\delta < t \leq 0\,,
     \end{cases}     
     \]
   whereas $\varphi = \varphi (x,v)$ is defined as 
	\[
	\begin{cases}
		\varphi \in C^\infty_c(B_{(\frac{\varrho +r}2)^{1+2s}}\times B_{\frac{\varrho +r}2}),\\
		0 \leq \varphi \leq 1~\text{and}~\phi \equiv 1~\text{on}~B_{{\varrho}^{1+2s}}\times B_{{\varrho}}\\
		|\nabla_v \varphi| \leq c/({r}-{\varrho})~\text{and}~|v\times \nabla_x\varphi| \leq c\langle v_{\rm o} \rangle/({r}-{\varrho})^{1+2s}.
	\end{cases}
	\] 
We observe that, by their definition, $\theta \in W^{1,2}([-r^{2s},0])$, whereas $\varphi$ is a smooth function. 

 Now, since we are working on a cylinder centered at the origin which is defined through euclidean open balls of suitable dimension (see \eqref{eq:cyl}), 
 we are allowed to introduce two symmetric standard mollifiers, the first one $\zeta_{h} = \zeta_{h} (t)$ in time supported in $(-h,h) \subset (t_1,t_2)$, and the second one $\gamma_{\ell}= \gamma_{\ell}(x)$ in space supported in $B_{\ell^{1+2s}}(0) 
\subset \Omega_x$, and 
	we define 
	\begin{align*}
		f_{h, \ell} (t,x,v) = \int_{\mathbb{R}}\int_{\mathbb{R}^n} \zeta_h(t-\tau) \gamma_\ell (x-\xi) f(\tau, \xi, v) \d\xi \d \tau.
	\end{align*}
	Note that in any domain $ \mathcal{U} \Subset (t_1, t_2) \times \mathbb{R}^{n}$ 
	and any $\mathcal{O} \subset \mathbb{R}^n$ on which 
	$f \in L^2(\mathcal{U}; H^s(\mathcal{O}))$, by the Lebesgue differentiation theorem for a.e. $t \in (t_1, t_2)$ we have that
	\begin{align*}
		\lim \limits_{h \to 0} \int_{\mathcal{U} \times \mathcal{O}} | f_{h,\ell}(t,x,v) - f_\ell(t,x,v) |^2 \, \dx \dv = 0 .
	\end{align*}
	Now, for any given~$\kappa \in \er$,  we denote $\phi (t,x,v)= \theta_\delta (t) \varphi(x,v)$ and introduce 
	$$
	 	\psi_{h,\ell} (t,x,v)= \left( (f_{h,\ell} - \kappa)_+ (t,x,v)  \phi(t,x,v) \right)_{h,\ell},
	$$
    	a test function which is $C^\infty_c (Q_r)$, and we test Definition~\ref{weak-sol-int} against it:
    	\begin{eqnarray*} \label{eq:subsol}
		0 &\geq& 
                  -  \int_{t_1}^{t_2}\int_{\Omega_x}
		\int_{\Omega_v} f (\partial_t+v\cdot\nabla_x)  \psi_{h,\ell} \dt\dx\dv  \\
        && 
                  + \int_{t_1}^{t_2}\int_{\Omega_x}\mathcal{E}(f,\psi_{h,\ell})\dt\dx 
                  =: J_{1, {h,\ell}} +J_{2, {h,\ell}}.	
      \end{eqnarray*}   
We begin estimating~$J_{1,{h,\ell}}$. Using the fact  that the convolution is symmetric and $f_\varepsilon$ is smooth with respect to $t$ and $x$, the term $J_{1, {h,\ell}}$ 
is equal to 
{
	\begin{align*}
	    J_{1, {h,\ell}} &=
	     -\int_{t_1}^{t_2}\int_{\Omega_x}
		\int_{\Omega_v} f_{h,\ell}(t,x,v) \phi(t,x,v) (\partial_t+v\cdot\nabla_x) (f_{h,\ell}-\kappa)_+ (t,x,v)  \dv\dx\dt\\ 
        &=
	     \int_{t_1}^{t_2}\int_{\Omega_x}
		\int_{\Omega_v} (\partial_t+v\cdot\nabla_x)f_{h,\ell}(t,x,v) \phi(t,x,v)  (f_{h,\ell}-\kappa)_+ (t,x,v)  \dv\dx\dt\\
	     &=  \frac12    \int_{t_1}^{t_2}\int_{\Omega_x}
		\int_{\Omega_v} (\partial_t+v\cdot\nabla_x) (f_{h,\ell}-\kappa)_+^2 (t,x,v) \phi(t,x,v) \dt\dx\dv \\
	     &= -\frac12 \int_{Q_r} (f_{h,\ell}-\kappa)_+^2 (t,x,v)(\partial_t+v\cdot \nabla_x) \phi(t,x,v) \dt\dx\dv.
	\end{align*}}		
Now, recalling the definition of $\varphi$ and $\theta_\delta$, 
and in particular since $\theta_\delta$ is a.e. differentiable and its derivative is equal to
 \[
 	\partial_t \theta_\delta=\frac{1}{\delta} \mathbf{1}_{(\tau_1-\delta,\tau_1]}- \frac{1}{\delta}\mathbf{1}_{(\tau_2,\tau_2+\delta]},
\] 
we infer
  \begin{eqnarray}
      -\frac12 \int_{Q_r} (f_{h,\ell}-\kappa)_+^2 (t,x,v)\varphi^2(x,v)\partial_t \theta_\delta(t) \dt\dx\dv & \geq & \frac{1}{2\delta}\int_{\tau_2}^{\tau_2+\delta}\int_{Q_\varrho^t} (f_{h,\ell}-\kappa)_+^2 (t,x,v)\dt\dx\dv\notag\\
      &&  - \frac{1}{2\delta}\int_{\tau_1-\delta}^{\tau_1}\int_{Q_r^t} (f_{h,\ell}-\kappa)_+^2 (t,x,v)\dt\dx\dv\notag
  \end{eqnarray}
 Then, putting the above computations into \eqref{eq:subsol}, and letting $\delta \to 0$ we recover the definition of derivative "in measure" and this yields
 \begin{eqnarray} \label{eq:j11}
		0 & \geq &
                 \frac{1}{2}\int_{B_{\varrho^{1+2s}\times B_\varrho}} (f_{h,\ell}-\kappa)_+^2 (\tau_2,x,v)\dt\dx\dv\\
      &&  - \frac{1}{2}\int_{B_{r^{1+2s}\times B_r}} (f_{h,\ell}-\kappa)_+^2 (\tau_1,x,v)\dt\dx\dv \notag\\
      &&
	    - \frac{c\langle v_{\rm o}\rangle}{(r-\varrho)^{1+2s}} \int_{B_{r^{1+2s}} \times B_r }(f_{h,\ell} -\kappa)_+^2(t,x,v) \, \dt\dx\dv \notag \\  
                  &&+ \int_{\tau_1}^{\tau_2}\int_{B_{r^{1+2s}}} \mathcal{E}(f,\left( (f_{h,\ell} - \kappa)_+ \varphi^2\right)_{h,\ell})\dt\dx\notag .	
      \end{eqnarray}   
Then, taking the averaged integral in $\tau_1 \in (-r^{2s}, 0]$ and taking $h,\ell \to 0$ 	

\begin{eqnarray} \label{eq:Cac1}
			&&  \int_{B_{\varrho^{1+2s}}\times B_\varrho} \omega^2(\tau_2,x,v)  \dx\dv  + \int_{-r^{2s}}^{\tau_2}\int _{B_{r^{1+2s}}} \mathcal{E}(f, \omega \varphi^2)\dt\dx\notag\\*
		  && \qquad\quad \leq  \frac{c\langle v_{\rm o}\rangle}{(r-\varrho)^{1+2s}} \int_{Q_r }\omega^2(t,x,v) \, \dt\dx\dv,
\end{eqnarray} 
where we have denoted with $\omega(t,x,v) := (f -\kappa)_+(t,x,v)$. Note that it would have been possible to consider a standard symmetric convolution 
defined according to the group (see for instance \cite{BLU07}), but in this case it would have been necessary to carry out the analysis on metric balls defined according to the
Carnot-Carathéodory distance of the group.

	Then, we estimate the second term on the left-hand side, that is the one related to the energy $ \mathcal{E}(f,\left( (f - \kappa)_+ \varphi^2\right) )$. Indeed, let us split the nonlocal energy as follows 
    \begin{eqnarray*}
    	&&\int _{B_r}
    	\int _{B_r}  \varPhi\left(\frac{f(t,x,v)-f(t,x,w)}{\snr{v-w}^{s}}\right) \notag\\*
        && \qquad \qquad\quad \times \frac{\big(\omega(t,x,v)\varphi^2(x,v) 
    	-\omega(t,x,w)\varphi^2(x,w)\big)}{\snr{v-w}^{n+s}}\dw \dv \nonumber
    	\\*[0.5ex]
    	&& + \, 2\int _{B_r} \int_{\er^n \setminus B_r} 
            \varPhi\left(\frac{f(t,x,v)-f(t,x,w)}{\snr{v-w}^{s}}\right) \frac{\omega(t,x,v)\varphi^2(x,v) }{\snr{v-w}^{n+s}} \dw\dv
    	\nonumber\\*[1.2ex]
    	&&=: J_{2,1} + J_{2,2},
    \end{eqnarray*}
    where we have recalled the definition of $\varphi$, which is only supported in $B_{r^{1+2s}} \times B_r$, and we  have considered that $\Phi$ is odd.
    
   We begin by estimating the argument of the term~$J_{2,1}$ by cases. With no loss of generality let us assume that $f(t,x,v) > f(t,x,w)$. If the opposite inequality holds true, then we exchange the roles of $v$ and $w$ by relying on the oddness of the nonlinearity $\varPhi$.
   
   Then, if $\omega(t,x,v)\varphi^2(x,v) 
    	- \omega(t,x,w)\varphi^2(x,w) \geq 0$,  by \eqref{eq:nonlinearity-ellipt-1} it holds
   $\varPhi(\tau)-\varPhi(\tau') \geq \Lambda^{-1}(\tau-\tau')$ for any $\tau \neq \tau'$, we obtain 
   \begin{eqnarray}\label{eq:cacc-lower-bound-1}
       && \varPhi\left(\frac{f(t,x,v)-f(t,x,w)}{\snr{v-w}^{s}}\right) \times \notag\\
       && \qquad \quad \times \frac{\Big(\omega(t,x,v)\varphi^2(x,v) 
    	- \omega(t,x,w)\varphi^2(x,w)\Big)}{\snr{v-w}^{n+s}}\nonumber\\
        && \quad \geq \Lambda^{-1} \Big(f(t,x,v)-f(t,x,w)\Big)\times \notag\\
        && \qquad\qquad \times \frac{\left(\omega(t,x,v)\varphi^2(x,v) 
    	- \omega(t,x,w)\varphi^2(x,w)\right)}{\snr{v-w}^{n+2s}}.
   \end{eqnarray}
   On  the other hand, if $\omega(t,x,v)\varphi^2(x,v) 
    	- \omega(t,x,w)\varphi^2(x,w) \leq 0$, by \eqref{eq:nonlinearity-ellipt-1} it holds
   $|\varPhi(\tau)| \geq \Lambda|\tau|$, so that 
   \begin{eqnarray}\label{eq:cacc-lower-bound-2}
       && \varPhi\left(\frac{f(t,x,v)-f(t,x,w)}{\snr{v-w}^{s}}\right) \times \notag\\
       && \qquad \quad \times \frac{\Big(\omega(t,x,v)\varphi^2(x,v) 
    	- \omega(t,x,w)\varphi^2(x,w)\Big)}{\snr{v-w}^{n+s}}\nonumber\\
        && \quad \geq \Lambda \Big(f(t,x,v)-f(t,x,w)\Big)\times \notag\\
        && \qquad\qquad \times\frac{\left(\omega(t,x,v)\varphi^2(x,v) 
    	- \omega(t,x,w)\varphi^2(x,w)\right)}{\snr{v-w}^{n+2s}}.
   \end{eqnarray}

Now we continue to split in different cases
{
   \begin{align*}
        &\left(f(t,x,v)-f(t,x,w)\right)\left(\omega(t,x,v)\varphi^2(x,v) 
    	- \omega(t,x,w)\varphi^2(x,w)\right) \\*[0.5ex] 
        & \quad =  \begin{cases} 
           \left(\omega(t,x,v)-\omega(t,x,w)\right) \notag\\
          \qquad \times \left(\omega(t,x,v)\varphi^2(x,v) 
    	- \omega(t,x,w)\varphi^2(x,w)\right) \, \, \, 
            &\text{if}\ f(t,x,v) > f(t,x,w) >\kappa,\\
           \left(f(t,x,v)-f(t,x,w)\right)\omega(t,x,v)\varphi^2(x,v) 
                \quad &\text{if} \ f(t,x,v) >\kappa \geq  f(t,x,w),\\
                0 \quad & \text{otherwise}
        \end{cases}\\
        &   \quad \geq \begin{cases} 
           \left(\omega(t,x,v)-\omega(t,x,w)\right) \notag\\
          \qquad \times \left(\omega(t,x,v)\varphi^2(x,v) 
    	- \omega(t,x,w)\varphi^2(x,w)\right) \, \, \, 
            &\text{if}\ f(t,x,v) > f(t,x,w) >\kappa,\\
          \omega^2(t,x,v)\varphi^2(x,v) 
                \quad &\text{if} \ f(t,x,v) >\kappa \geq  f(t,x,w),\\
                0 \quad & \text{otherwise}
        \end{cases}\\
         &   \quad \geq \big((\omega\varphi)(t,x,v)-(\omega\varphi)(t,x,w)\big) ^2 -\omega(t,x,v)\omega(t,x,w)\big(\varphi(x,v)-\varphi(x,w)\big)^2.
    \end{align*}}
 Combining the above estimates into $J_{2,1}$ yields 
      \begin{eqnarray*}
    J_{2,1}
    & \geq &  c  [\omega\varphi]^2_{H^s(B_r)}  \\
    &-& c \iint_{B_{r}\times B_{r}}\tfrac{\max\{ \omega(t,x,v),\omega(t,x,w) \}^2\snr{\F(x,v)-\F(x,w)}^2}{\snr{v-w}^{n+2s}} \dv \dw 
    \end{eqnarray*}
    for some $c \equiv c(\Lambda)>0$. Now, we apply further estimate the last term. By symmetry of the Gagliardo kernel we can assume with no loss of generality that $\omega(t,x,v) \geq \omega(t,x,w)$, up to exchanging the roles of $v$ and $w$. Hence,
    \begin{eqnarray} \label{eq:max}
        &&  \iint_{B_{r}\times B_{r}}\frac{\max\{ \omega(t,x,v),\omega(t,x,w) \}^2\snr{\F(x,v)-\F(x,w)}^2}{\snr{v-w}^{n+2s}} \dv \dw \notag \\
        && \quad \leq \int_{B_r} \omega^2(t,x,v)\left(\int_{B_{2r}(v)}\frac{\|\nabla_v\varphi\|^2_{L^\infty(Q^t_r)}\dw}{\snr{v-w}^{n-2(1-s)}} \right) \dv \notag\\
        && \quad \leq \frac{c\,r^{2(1-s)}}{(r-\rr)^2}\int_{B_r}\omega^2(t,x,v)\dv. 
    \end{eqnarray}
   All in all, combining the above estimates yields 
   \begin{eqnarray}\label{eq:J11-final}
    J_{2,1}
    & \geq &  c [\omega\varphi]^2_{H^s(B_r)}  - \frac{c}{(r-\rr)^2}\int_{B_r}\omega^2(t,x,v)\dv .
    \end{eqnarray}

    Now, we deal with the nonlocal term in~$J_{2,2}$. Firstly, we observe that when let us assume $f(t,x,v) > f(t,x,w)$, then by the second estimate in
   \eqref{eq:nonlinearity-ellipt-1} we get
    \begin{eqnarray*}
        &&\varPhi\left(\frac{f(t,x,v)-f(t,x,w)}{\snr{v-w}^{s}}\right) \frac{\omega(t,x,v)\varphi^2(x,v) }{\snr{v-w}^{n+s}} \\
      %  &&\qquad =
      %  \varPhi\left(\frac{f(t,x,v)-f(t,x,w)}{\snr{v-w}^{s}}\right) \frac{f(t,x,v) - f(t,x,w)}{\snr{v-w}^s} \frac{ \omega(t,x,v)\varphi^2(x,v) }{(f(t,x,v) - f(t,x,w))\snr{v-w}^{n}} \\
        &&\qquad \geq \Lambda^{-1} \frac{f(t,x,v) - f(t,x,w)}{\snr{v-w}^{n+2s}} \omega(t,x,v)\varphi^2(x,v) \\
        &&\qquad \geq \Lambda^{-1} \frac{\kappa - f(t,x,w)}{\snr{v-w}^{n+2s}} \omega(t,x,v)\varphi^2(x,v) \\
         &&\qquad \geq -\Lambda^{-1} \frac{\omega(t,x,w)\omega(t,x,v)\varphi^2(x,v) }{\snr{v-w}^{n+2s}}.
    \end{eqnarray*}
    Note that the above chain of inequalities holds trivially when $f(t,x,v) = f(t,x,w)$. Finally, when $f(t,x,v) < f(t,x,w)$ we get an analogous estimate by considering the first estimate in Assumption \ref{assumption} as follows 
    \begin{eqnarray*}
        &&\varPhi\left(\frac{f(t,x,v)-f(t,x,w)}{\snr{v-w}^{s}}\right) \frac{\omega(t,x,v)\varphi^2(x,v) }{\snr{v-w}^{n+s}} \\
        &&\qquad=\varPhi\left(\frac{f(t,x,w)-f(t,x,v)}{\snr{v-w}^{s}}\right) \left( - \frac{\omega(t,x,v)\varphi^2(x,v) }{\snr{v-w}^{n+s}}\right) \\
        &&\qquad \geq \Lambda \frac{f(t,x,w) - f(t,x,v)}{\snr{v-w}^{n+2s}} ( - \omega(t,x,v)\varphi^2(x,v))  \\
          &&\qquad \geq \Lambda \frac{f(t,x,w) - \kappa}{\snr{v-w}^{n+2s}} ( - \omega(t,x,v)\varphi^2(x,v))  \\
        &&\qquad \geq - \Lambda \frac{\omega(t,x,w)\omega(t,x,v)\varphi^2(x,v) }{\snr{v-w}^{n+2s}}.
    \end{eqnarray*}
   Then~$J_{2,2}$ can be estimated as follows:
   \begin{equation}\label{eq:J21}
   	J_{2,2} 	 \geq -c \int_{B_r} \int_{\er^n \setminus B_r}
   	\frac{\varphi^2(x,v)\omega(t,x,v)\omega(t,x,w)}{\snr{v-w}^{n+2s}} \dw\dv. 
   \end{equation}

   Now, combining \eqref{eq:J11-final} and \eqref{eq:J21} with \eqref{eq:Cac1}
   \begin{eqnarray*}
   	&&  \int_{B_{\varrho^{1+2s}}\times B_\varrho} \omega^2(\tau_2,x,v)  \dx\dv  + \int_{-r^{2s}}^{\tau_2}\int_{B_{r^{1+2s}}}[\omega\varphi]^2_{H^s(B_r)}\dt\dx\notag\\*
		  && \quad \leq  \frac{c\langle v_{\rm o}\rangle}{(r-\varrho)^{2(1+s)}} \int_{Q_r }\omega^2(t,x,v) \, \dt\dx\dv\\
          && \qquad + c \int_{Q_r} \int_{\er^n \setminus B_r}
   	\frac{\varphi^2(x,v)\omega(t,x,v)\omega(t,x,w)}{\snr{v-w}^{n+2s}} \dw\dv\dx\dt\notag\,,
  \end{eqnarray*}
   where the constant $c>0$ depends only on $\Lambda$, $n$ and $s$.
   Now, with some standard manipulations as in \cite[Lemma 3.1]{APP25}, we arrive at
   \begin{eqnarray}\label{eq:Cacc-final-2}
   	&&  \sup_{t \in [-\rr^{2s},0]} \int_{B_{\varrho^{1+2s}}\times B_\varrho} \omega^2(t,x,v)  \dx\dv  + \int_{U_\rr}[\omega]^2_{H^s(B_\rr)}\dt\dx\notag\\*
		  && \quad \leq  \frac{c\langle v_{\rm o}\rangle}{(r-\varrho)^{2(1+s)}} \int_{Q_r}\omega^2(t,x,v) \, \dt\dx\dv\\
          && \qquad + c \int_{Q_r} \int_{\er^n \setminus B_r}
   	\frac{\varphi^2(x,v)\omega(t,x,v)\omega(t,x,w)}{\snr{v-w}^{n+2s}} \dw\dt\dx\dv\notag\,.
  \end{eqnarray}
   Furthermore, we estimate from above the nonlocal tail on the right-hand side of \eqref{eq:Cacc-final-2} by 	 applying H\"older's Inequality  with~$\big(\sigma,\frac{\sigma}{\sigma-1}\big)$, with~$\sigma > 2$. In this way we obtain
{	\begin{align*}
		&\int_{Q_{{r}}}\int_{\ern \setminus B_r} \frac{\varphi^2(x,v)\omega(t,x,v)\omega(t,x,w)}{\snr{v-w}^{n+2s}} \dw\dt\dx\dv\nonumber\\
		& \leq \int_{Q_{{r}}\cap \text{supp}(\varphi)} \varphi^2(x,v)\omega(t,x,v) \left(\int_{\ern \setminus B_r}
		\frac{\omega(t,x,w)}{\snr{v-w}^{n+2s}}\dw\right)\dt\dx\dv\nonumber\\
		& \leq   \left(\int_{{Q}_{{r}}\cap \{f>\kappa\}}\omega^\frac{\sigma}{\sigma-1} \dt\dx\dv\right)^\frac{\sigma-1}{\sigma} \times \nonumber\\*[1ex]
		&\hspace{2cm}\times \left[\int_{{Q}_{{r}}\cap \, 
		 \text{supp}(\varphi)}\left(\int_{\ern \setminus B_r}\frac{ \omega(t,x,w)}{\snr{v-w}^{n+2s}}\dw\right)^\sigma \dt \dx \dv\right]^\frac{1}{\sigma}\nonumber\\*[1ex]
		& \leq   \frac{c\,|Q_r \cap \{f > \kappa\}|^{\frac 12-\frac{1}{\sigma}}}{({r}-{\varrho})^{n+2s}}\|\omega\|_{L^2(Q_r)} \|\tail(\omega;B_{{r}})\|_{L^\sigma(U_r)},
	\end{align*}
	where in the last display we applied H\"older's Inequality once again with $\big(\frac{2(\sigma-1)}{\sigma}, \frac{2(\sigma-1)}{\sigma-2}\big)$, noting that $\sigma>2$ implies that $\sigma/(\sigma-1)<2$} and we centered the Gagliardo kernel since for any~$v \in B_{{r}} \cap \textup{supp}(\varphi) \subset B_{({r}+{\varrho})/2}$ and any~$w \in \ern \setminus B_r$, it  holds
	\[
	\frac{|w|}{|v-w|}\, \leq \,1 + \frac{|v|}{||w|-|v||} \,\leq\, 1 + \frac{{r}+{\varrho}}{{r}-{\varrho}}\, =\, \frac{c\,{r}}{{r}-\varrho}.
	\]
  Then, combining the above estimates with \eqref{eq:Cacc-final-2} yields the desired result.
  \end{proof}

  \begin{proof}[\bf Proof of Theorem~\ref{loc bdd}]
  The proof follows with an analogous procedure as in \cite[Theorem 1.1]{APP25}. For the sake of the reader we just give a sketch 
  here as well. Note now, that by Assumption \ref{assumption} it follows that
  \begin{equation} \label{eq:riscr1}
       \Lc f(t,x,v) := {\rm P.V.}\int_{\ern}\big(f(t,x,v)-f(t,x,w)\big)K(t,x,v,w) \dw
  \end{equation}
  where for every $v \ne w$ the symmetric kernel $K$ is defined as
  \begin{equation} \label{eq:riscr2}
  K(t,x,v,w) = \varPhi\left(\frac{f(t,x,v)-f(t,x,w)}{|v-w|^s}\right)\frac{|v-w|^{-n-s}}{f(t,x,v)-f(t,x,w)},
  \end{equation}
  and satisfies 
  \begin{equation} \label{eq:riscr3}
  \frac{\Lambda^{-1}}{|v-w|^{n+2s}} \leq K(t,x,v,w) \leq \frac{\Lambda}{|v-w|^{n+2s}}\,,
  \end{equation}
  see \cite[Remark 4.1]{DKLN25a}. 
  Now, the proof proceed in two steps.

  \subsection*{Step 1: The gain of integrability}
   By performing the same argument as in \cite{APP25} -- see also \cite{Hou25} for a related approach in the local case -- we can rely on the higher integrability estimates achievable via the fundamental solution of the fractional Kolmogorov equation, in turn obtaining the following Sobolev type inequality (see \cite[Theorem 1.4]{APP25})
\begin{eqnarray}  \label{eq:gain}
	&&  \hspace{-5mm}{(r-\varrho)^{n+2s}}\| (f-\kappa)_+\|_{L^q ({Q}_{\varrho}(z_{\rm o}))} \notag\\*[0.5ex]
	&& \quad \leq   c\,\langle v_{\rm o}\rangle  \|(f-\kappa)_+ \|_{L^2({U}_{r}(t_{\rm o},x_{\rm o}); H^s(B_r(v_{\rm o})))} \\*
	&& \qquad 	+\,c\,|{Q}_{r}(z_{\rm o}) \cap \{f >\kappa\}|^{\frac{1}{2}+\frac{s}{N_{s}} -\frac{1}{\sigma}}\|\textup{Tail}((f-\kappa)_+;B_{r}(v_{\rm o}))\|_{L^{\sigma}(U_{r}(t_{\rm o},x_{\rm o}))} \notag\,,
\end{eqnarray} 
for any~$\kappa \in \er$, any~$\varrho \in (0,r)$ where the constants~$c\equiv c(n,s,\Lambda,p)>0$ and the exponent $(\sigma,q)$ satisfies
$$
	2 \leq q \leq \frac{2N_{s}}{N_{s}-2s}\,, \qquad \text{and} \qquad \sigma > \frac{N_{s}}{2s}\,.
$$

\subsection*{Step 2: De Giorgi iteration}
The second step of the proof is based on a classical De Giorgi argument. Indeed, let us translate the problem considering~$\tilde{f}(z):=f(z_{\rm o}\circ z)$, 
so that the center of the cylinder is the origin. Hence, by~\cite[Lemma~5.1]{Sto19} we have that~$\tilde{f}$ solves 
\[
(\partial_t + v\cdot \nabla_x)\tilde{f}= \tilde{\Lc} \tilde{f} \quad \text{in}~ \tilde{\Omega}:= z_{\rm o}^{-1}\circ\Omega\,,
\]
where~$\tilde{\Lc}$ is an integro-differential operator whose kernel satisfies the same ellipticity condition as in~\eqref{eq:riscr3}.
Then, for any~$j \in \mathbb{N}$, define~
$$
r_j :=  \frac{1}{2}(1+2^{-j})r \quad \text{and}\quad \kappa_j:=(1-2^{-j})\kappa\,,
$$
where~$\kappa>0$ will be fixed later on. 
 Now, we apply estimate~\eqref{eq:gain} to~$(\tilde{f}-\kappa_{j+1})_+$, with~$q =\frac{2N_{s}}{N_{s}-2s}$ and with radii~$r_{j+1}$ and~$r_j$.

Let us first note that by the Caccioppoli estimate \eqref{caccioppoli-def} we can estimate the $L^2({U}_{r}(t_{\rm o},x_{\rm o}); H^s(B_r(v_{\rm o})))$-norm of $(f-\kappa_{j+1})_+$ on the right-hand side of \eqref{eq:gain}.

Indeed, first, by Chebychev's Inequality we have
\begin{eqnarray}\label{eq:cheb-bdd}
   \frac{|Q_{r_j} \cap \{\tilde{f}> \kappa_{j+1}\}|}{|Q_{r_j}|} & \leq &   \,\frac{c\,2^{2j}}{\kappa^2}\mean{Q_{r_j}}(\tilde{f} -\kappa_j)_+^2\dt\dx\dv.
\end{eqnarray}
Indeed, first of all, let us note that choosing
\begin{equation}\label{eq:k0}
	\kappa \,\geq\, \delta  \nra{{\tail((\tilde{f})_+;B_{r/2})}}_{L^\sigma(U_r)}\qquad \text{for}~\delta \in (0,1]\,,
\end{equation}
where the notation of the right-hand side is defined in \eqref{nra},
yields 
\begin{eqnarray*}
&& \frac{|Q_{r_j}\cap \{\tilde{f} > \kappa_j\}|^{\frac12-\frac{1}{\sigma}}}{|Q_{r_j}|}
        \|\tail((\tilde{f}-\kappa_j)_+;B_{{r_j}}) \|_{L^\sigma ( U_{r_j}) }\left( \int_{Q_{r_j}} (\tilde{f}-\kappa)_+^2\dt \dx \dv \right)^\frac12\notag\\
&& \quad \leq |B_{r_j}|^{-\frac{1}{\sigma}}\left(\frac{|Q_{r_j}\cap \{\tilde{f} > \kappa_j\}|}{|Q_{r_j}|}\right)^{\frac12-\frac{1}{\sigma}}\nra{\tail((\tilde{f}-\kappa_j)_+;B_{{r_j}})}_{L^\sigma ( U_{r_j}) }\times\\ 
&& \qquad\qquad\qquad\quad \times \left(\mean{Q_{r_j}} (\tilde{f}-\kappa)_+^2\dt \dx \dv \right)^\frac12\\
&& \quad \leq  c2^{j}\kappa|B_{r_j}|^{-\frac{1}{\sigma}}\left(\mean{Q_{r_j}} \frac{(\tilde{f}-\kappa)_+^2}{\kappa^2}\dt \dx \dv\right)^{1-\frac{1}{\sigma}}\nra{\tail(\tilde{f})_+;B_{\frac{r}{2}})}_{L^\sigma ( U_{r}) }\notag\\
&& \quad \leq c\,r^{- \frac{n}{\sigma}}2^{j}\left(\frac{\kappa}{\delta}\right)^2\left(\mean{Q_{r_j}} \frac{(\tilde{f}-\kappa)_+^2}{\kappa^2}\dt \dx \dv\right)^{1-\frac{1}{\sigma}}.
\end{eqnarray*}

Thus, by the Caccioppoli estimate \eqref{caccioppoli-def} we obtain
\begin{eqnarray}\label{eq:nonloc-en-bdd}
 && \mean{Q_{r_{j+1}}}\int_{B_{r_{j+1}}}\frac{\snr{(\tilde{f}-\kappa)_+(t,x,v)
 -(\tilde{f}-\kappa)_+(t,x,w)}^2}{\snr{v-w}^{n+2s}}\dw \dt \dx \dv \notag\\*[0.5ex] 
   && \quad \  \leq\ {c\, r^{- \frac{n}{\sigma}-  2(n+2s)}  2^{2j(1+s)+ jN_{s}}\langle{v_{\rm o}\rangle^2}\nra{{(\tilde{f} -\kappa_j)_+}}_{L^2(Q_{r_j})}^2
  } \\*[0.5ex]
   && \qquad \quad + c r^{- \frac{n}{\sigma}-  2(n+2s)} 2^{4j(n+s)+jN_{s}}\left(\frac{\kappa}{\delta}\right)^2\nra{{(\tilde{f} -\kappa_j)_+/\kappa}}_{L^2(Q_{r_j})}^{2(1-\frac{1}{\sigma})}  \,,\notag
\end{eqnarray}
where we used also that $|Q_{r_{j}}|/|Q_{r_{j+1}}| \lesssim 2^{jN_{s}}$.

Also, from \eqref{eq:cheb-bdd} and the choice of $\kappa$ in \eqref{eq:k0} we obtain
\begin{eqnarray}\label{eq:tail-est}
		&&  \|{\textup{Tail}((\tilde{f} -\kappa_{j+1})_+;B_{r_j})}\|_{L^\sigma(U_{r_j})}^2
        \frac{|Q_{r_j}\cap \{\tilde{f}> \kappa_{j+1}\}|^{1+\frac{2s}{N_{s}}-\frac{2}{\sigma}}}{\snr{Q_{r_j}}}\notag\\
        	&& \quad \leq  |Q_{r_j}|^{\frac{2s}{N_s}}|B_{r_j}|^{-\frac{2}{\sigma}}\nra{\textup{Tail}((\tilde{f} -\kappa_{j+1})_+;B_{r_j})}_{L^\sigma(U_{r_j})}^2\times \notag\\
            && \qquad\qquad\qquad\quad \qquad \times
        \left(\frac{|Q_{r_j}\cap \{\tilde{f}> \kappa_{j+1}\}|}{\snr{Q_{r_j}}}\right)^{1+\frac{2s}{N_{s}}-\frac{2}{\sigma}}\notag\\
        && \quad \leq  c\,r^{-\frac{2n}{\sigma}} 2^{j(2+\frac{4s}{N_s} -\frac{4}{\sigma})}|Q_1|^{\frac{2s}{N_s}}\nra{\textup{Tail}((\tilde{f})_+;B_{\frac r2})}_{L^\sigma(U_r)}^2\times \notag\\
        && \qquad\qquad\qquad\quad \times
        \left(\mean{Q_{r_j}}\frac{(\tilde{f} -\kappa_j)_+^2}{\kappa^2}\dt\dx\dv\right)^{1+\frac{2s}{N_{s}}-\frac{2}{\sigma}}\notag\\
         && \quad \leq  c\,r^{-\frac{2n}{\sigma}} 2^{j(2+\frac{4s}{N_s} -\frac{4}{\sigma})}\left(\frac{\kappa}{\delta}\right)^2
     \nra{{(\tilde{f} -\kappa_j)_+/\kappa}}_{L^2(Q_{r_j})}^{2(1+\frac{2s}{N_{s}}-\frac{2}{\sigma})}.
	\end{eqnarray}

Then, combining~\eqref{eq:tail-est} and~\eqref{eq:nonloc-en-bdd} together with~\eqref{eq:gain} yields 
		\begin{eqnarray*}
			&& \hspace{-6mm}\nra{ {(\tilde{f} -\kappa_{j+1})_+}}_{L^2(Q_{r_{j+1}})}^2  \notag\\
			&& \quad \leq    \frac{\nra{(\tilde{f}-\kappa_{j+1})_+}_{L^q(Q_{r_j})}^2}{|Q_{r_j}|}
         |Q_{r_{j+1}} \cap \{\tilde{f}> \kappa_{j+1}\}|^{\frac{2s}{N_s}} \notag\\*
			&& \quad \leq c_*\,b^j \kappa^2\Biggl[Y_j + Y_j^{1-\frac{1}{\sigma}} + Y_j^{1+\frac{2s}{N_{s}}-\frac{2}{\sigma}}\Biggr] \left(\frac{|Q_{r_{j+1}} \cap \{\tilde{f}> \kappa_{j+1}\}|}{\snr{Q_{r_{j+1}}}}\right)^{\frac{2s}{N_{s}}}\,,
		\end{eqnarray*}
		for
		\begin{eqnarray*}
		&&    b \equiv b(n,s)>1,
		\qquad   c_* := \big({c\,\delta^{-1}r^{-\frac{n}{\sigma} -  2(n+2s)}}\big)^2\langle v_{\rm o}\rangle^3>0,\\*[0.5ex]
         \text{and} \quad && Y_j:= \frac{\nra{{(\tilde{f} -\kappa_j)_+}}_{L^2(Q_{r_j})}^2}{\kappa^2}.
		\end{eqnarray*}
	Thus, by applying once again Chebychev's Inequality \eqref{eq:cheb-bdd}, up to eventually relabeling $b$ and $c$ we get
	\begin{equation}\label{eq:bdd-iter}
		Y_{j+1}\, \leq c_*\,b^j\left(Y_j^{1+\frac{2s}{N_{s}}} + Y_j^{1+2(\frac{2s}{N_{s}}-\frac{1}{\sigma})}  + Y_j^{1+\frac{2s}{N_{s}}-\frac{1}{\sigma}}\right).
	\end{equation}
	Note  that~${N_{s}}/({2s})<\sigma$ implies that~$\frac{2s}{N_{s}}> \frac{1}{\sigma}$.

	Hence, up to choosing~$\kappa$ such that
	\begin{equation}\label{eq:k0-2}
		\kappa  \geq \nra{ {\tilde{f}_+}}_{L^2(Q_r)} \,,
	\end{equation} 
	we can rewrite~\eqref{eq:bdd-iter} as follows
	\[
	Y_{j+1}\, \leq\, c_* b^j Y_{j}^{1+{\alpha}}\,,
	\]
	for some positive~$\alpha \equiv\alpha(n,s,\sigma) := \frac{2s}{N_{s}}-\frac 1\sigma>0$ and~$b>1$. Then, up to choosing (upon translating and dilating back)
	\begin{align*}
		& \kappa  := 
		b^{\frac{1}{2\alpha^2}}c^\frac{1}{\alpha}\frac{\langle v_{\rm o}\rangle^{\frac{3}{2\alpha}}}{\left( r^{\frac{n}{\sigma}+2(n+2s)}
		\delta \right)^{\frac1\alpha} }\nra{f_+}_{L^2(Q_r(z_{\rm o})} + \ \delta  \, 
		\nra{ {\tail(f_+;B_\frac{r}{2}(v_{\rm o}))}}_{L^\sigma(U_{r}(t_{\rm o},x_{\rm o}))}\,,
	\end{align*}
	in clear accordance with~\eqref{eq:k0} and~\eqref{eq:k0-2}, the iteration argument of 
	\cite[Lemma 2.6]{DKP14} yields that~$Y_j \to 0$ as~$j \to \infty$, which gives the desired result.
\end{proof}

\section{Further comments on general nonlinear nonlocal diffusions}
\label{sec:nonlinear}
In this section, we give a brief overview of the particular case when $\varPhi(\tau)=|\tau|^{p-2}\tau$,
and in particular we highlight the major challenges one has to face 
to deal with nonlocal kinetic $p$-Laplace equations in contrast to the elliptic case~\cite{DKP14,DKP16},
or the parabolic one~\cite{APT22,Liao22,Tav24}. Our equation reads as follow
\begin{equation}\label{pbm4}
		(\partial_t +v \cdot \nabla_x) f = \mathcal{L} f \qquad \text{for}~(t,x,v) \in \er \times \ern \times \ern\,,
	\end{equation}
where the diffusion term~$\mathcal{L}$ is an integro-differential operator of differentiability order~$s \in (0,1)$ and summability order~$p\in (1,\infty)$ given by {
\begin{align}\label{diff4}
			\Lc f(t,x,v) := {\rm P.~\!V.}\int _{\ern} 
			&\snr{f(t,x,v)-f(t,x,w)}^{p-2}  (f(t,x,w)-f(t,x,v))K(t,x,v,w) \dw,
\end{align}}
where $K$ is a symmetric measurable kernel such that
	\begin{equation}\label{main krn}
		\Lambda^{-1}\snr{v-w}^{-n-sp} \leq \K(t,x,v,w) \leq \Lambda \snr{v-w}^{-n-sp} \quad \text{for a.~\!e.}~v,w \in \ern,
  	\end{equation}  
 for a.~\!e.~$(t,x)\in \er^{1+n}$ and for a positive constant~$\Lambda> 0$.
 	
 \vspace{2mm}
 
 As a prototype for Equation~\eqref{pbm4}, even though in this scenario the difficulties arising when dealing with only measurable coefficients vanishes, one can consider the simpler case when the involved kernel~$\K$ does coincide with the classical Gagliardo kernel in velocity, i.~\!e.~$ \K(t,x,v,w) \equiv \snr{v-w}^{-n-sp}$. In this setting, Equation~\eqref{pbm4} does reduce to
$$
	(\partial_t +v \cdot \nabla_x) f +  (-\Delta_v)_p^sf=0,
$$
where~$(-\Delta_v)_p^s$ is the classical $(s,p)$-Laplacian with respect to the $v$-variable. 

\vspace{2mm}

Employing the techniques already proposed in this work, see Section \ref{sec_cc}, 
we are able to prove a Caccioppoli inequality for weak solutions of the equation above. 
For this, one has to introduce the proper notion of weak solution and the correct geometry to deal with 
the new nonlinear fractional setting. 

 In a similar fashion of what already done in Section \ref{sec:preliminaries}, in this case we endow~$\er^{1+2n}=\er\times\er^n\times \er^n$ with the same group law $\circ$, whereas for any~$p>1$ and for any $r>0$, we consider the usual fractional nonlinear kinetic scaling~$\delta_r:  \er^{1+2n} \mapsto \er^{1+2n}$ defined by
\begin{equation}\label{eq:dil2}
	\delta_r(t,x,v):=(r^{sp}t,\, r^{1+sp}x,\, r v). 
\end{equation}
Then for any~$r>0$, the {\it slanted} cylinder~${Q}_r(z_{\rm o})$ is defined as follows,
\begin{eqnarray*}
	&& {Q}_r(z_{\rm o})\!\!\! \ :=
	\big \{ z:=(t,x,v) \in \er^{1+2n}: \, -r^{2p} < t -  t_{\rm o} \leq  0 , \\* 
&&\hspace{3,5cm}   \ | x - x_{\rm o} - (t - t_{\rm o}) v_{\rm o}| <r^{1 + 2p},\ \snr{ v- v_{\rm o} } < r \big\}, \notag
\end{eqnarray*}
and the {\it homogeneous dimension} $N_{sp}$ related to~\eqref{eq:dil2} is defined as
\[	N_{sp}:=  n(2+sp)+sp.
\]

Furthermore, for~$p \in (1,\infty)$,~$s \in (0,1)$ and any~$\mathcal{O}\subseteq \ern$, we denote with~$W^{s,p}(\mathcal{O})$ the fractional Sobolev space
	\[
	W^{s,p}(\mathcal{O}) := \big\{ f  \in L^p(\mathcal{O}): [f]_{W^{s,p}(\mathcal{O})} < +\infty\big\},
	\]
	where the fractional seminorm~$[f]_{W^{s,p}(\mathcal{O})}$ is the usual one via Gagliardo kernels, 
	\[
	[f]_{W^{s,p}(\mathcal{O})}:= \left(\iint_{\mathcal{O}\times \mathcal{O}}\frac{\snr{f (v)-f (w)}^p}{\snr{v-w}^{n+sp}}\dv \dw\right)^\frac{1}{p}.
	\]
	We endow~$W^{s,p}(\mathcal{O})$ with the following norm
	\[
		\|f\|_{W^{s,p}(\mathcal{O})} := \|f\|_{L^p(\mathcal{O})} + [f ]_{W^{s,p}(\mathcal{O})}.
   	 \]
   	 A function~$ f$ belongs to $W_{\rm loc}^{s,p}(\mathcal{O})$ if~$f\in W^{s,p}(\mathcal{O}')$ whenever $\mathcal{O}' \subset\subset \mathcal{O}$.
	In a similar fashion, we denote with~$W^{s,p}_0(\mathcal{O})$ the closure of~$C^\infty_0(\mathcal{O})$ with respect to~$\norma{\times}_{W^{s,p}(\mathcal{O})}$.
	Lastly, we recall that as a definition of tail we consider the nonlinear version of the one suggested in~\eqref{eq:tail}
    $$
     \tail(f;B_r(v_{\rm o})):= r^{sp}\int_{\ern \setminus B_r(v_{\rm o})} \frac{\snr{f(t,x,v)}^{p-1}}{\snr{v_{\rm o}-v}^{n+sp}}\dv\,.
    $$
	and we consider the corresponding {\it tail space}
\[
L^{p-1}_{sp}(\ern):= \left\{g \in L^{p-1}_{\textrm{loc}}(\ern)\, : \,  \|g\|_{L^{p-1}_{sp}(\ern)}:= \int_{\ern}\frac{|g(v)|^{p-1}}{(1+\snr{v})^{n+sp}}\dv < \infty \right\},
\]
as firstly defined in~\cite{KKP16}; see Section~2 there for related properties.

\vspace{2mm}

Then, following the steps of the classical approach for instance proposed in \cite{Strom19}, 
given~$\Omega:= (t_1,t_2)\times \Omega_x \times \Omega_v \subset \er^{1+2n}$  
we denote by~$\Wc$ the natural functions space to which weak solutions to~\eqref{pbm} of our interest belong to, and we defined it as
\begin{multline*}
	\W  :=   \Big\{ f  \in L^p  ((t_1,t_2)\times \Omega_x ;W^{s,p}(\Omega_v)) \cap L^{p-1}((t_1,t_2)\times \Omega_x; L^{p-1}_{sp}(\ern))\\
	 \, 
	: (\partial_t + v \times\nabla_x)f  \in L^{p'}((t_1,t_2)\times \Omega_x ; ( W^{s,p}(\ern) )^* ) \Big \},
\end{multline*}
where $( W^{s,p}(\ern) )^*$ is the dual space of $W^{s,p}(\ern)$, and if~$ p \in (1, + \infty)$ we denote~$p':=p/(p-1)$ as its conjugate exponent.

Furthermore, we denote by~$\Ec$ the nonlocal energy associated with our diffusion term~$\mathcal{L}$ in~\eqref{diff4}; that is {
\begin{align*}
	\Ec(f,\varphi) := \iint_{\ern\times\ern} &\snr{f(v)-f(w)}^{p-2}  (f(v)-f(w))
	(\varphi(v)-\varphi(w))K(v,w)\, \dv\dw\,,
\end{align*}}
for any test function~$\varphi$ smooth enough.
We are now in a position to recall the definition of weak sub- and supersolution.
\begin{defn}\label{weak-sol-int-2}
	A function~$f \in \Wc$ is a \,{\rm weak subsolution} {\rm (}resp.,~{\rm supersolution}{\rm)} to~\eqref{pbm4} in~$\Omega$ if 
	\begin{eqnarray*}
		&&\int_{t_1}^{t_2}\int_{\Omega_x}  \Ec(f,\varphi) \dt\dx \\
		&&\quad - \,\int_{t_1}^{t_2}\int_{\Omega_x}
		\int_{\Omega_v}  \,f \,  (\partial_t+v\cdot\nabla_x)  \varphi  \dt\dx\dv
		\, \leq\, 0\ \ \big(\ge\,0, \, \text{resp.}\big),
	\end{eqnarray*}
	for any nonnegative~$\varphi \in L^p  ((t_1,t_2)\times \Omega_x ;W^{s,p}(\ern))$ 
	and $(\partial_t+v\times\nabla_x) \varphi \in L^{p'} ((t_1,t_2)\times \Omega_x ; (W^{s,p}(\ern))^*)$,
	such that~${\rm supp} \,\varphi \subset\subset (t_1, t_2) \times \Omega_x \times \Omega_v$.
	%and where $\langle \times | \times \rangle$ is the usual duality pairing in velocity
    \\
	A function~$f \in \Wc^p$ is a {\rm weak solution} to~\eqref{pbm4} if it is both a weak sub- and supersolution.
\end{defn}
Note that, when a general $p>1$ is considered, we cannot get rid of 
the requirement ${\rm supp} \,\varphi \Subset \Omega$, because it ensures the validity 
of the integration by parts formula we will later on employ.

\begin{lemma}\label{thm ccpp}
   	Let~$p \in (1, + \infty)$, $s \in (0,1)$ and $\Omega\subset \er^{1+2n}$ be a bounded domain.
   Let~$f $ be a weak subsolution to~\eqref{pbm4} in~$\Omega$ according to {\rm Definition~\ref{weak-sol-int-2}}. For any $r \in (0,1)$ such that~$Q_{r}(z_{\rm o}) \subset \Omega$ the following estimate holds true for any $\kappa \in \er$, $\sigma>\frac{p}{p-1}$ and $\rr \in (0,r)$
   \begin{eqnarray*} \nonumber
   	 	&& \sup_{t \in (-\rr^{sp}+t_{\rm o},t_{\rm o}]}
		\int_{Q_\rr^t(z_{\rm o})} (f-\kappa)^2_+\dx\dv  +  \int_{U_\rr(t_{\rm o},x_{\rm o})} [(f-\kappa)_+]^2_{W^{s,p}(B_\rr(v_{\rm o}))} \dt \dx \nonumber \\*[0.5ex]
   		&& \quad\leq   \frac{c\langle v_{\rm o}\rangle} 
              {(r-\rr)^{p(1+s)}}   \int_{Q_r(z_{\rm o})} (f-\kappa)_+^p\dt \dx \dv \\
        &&\quad +\frac{c\,|Q_r(z_{\rm o})\cap \{f > \kappa\}|^{1-\frac1p-\frac{1}{\sigma}}}{({r}-{\varrho})^{2(n+sp)}}\|(f-\kappa)_+\|_{L^p(Q_r(z_{\rm o}))}
        \|\tail((f-\kappa)_+;B_{{r}}(v_{\rm o})) \|_{L^\sigma\left( U_r(t_{\rm o},x_{\rm o})\right) }\,, \nonumber
   	\end{eqnarray*}
	where the constant~$c>0$ depends only on $p$ and on the kernel constant $\Lambda$.
\end{lemma}

\begin{proof}
With no loss of generality let us assume $z_{\rm o}=0$.
      Let~$Q_r \equiv Q_r({0}) \Subset \Omega$ and let~$f$ be a weak subsolution to~\eqref{pbm} according to Definition~\ref{weak-sol-int}.  
      
      The proof is carried out as in Lemma \ref{thm ccp},
      %the only difference here is that the mollifier needs to be constructed according to Lie group given by $\mathbb{R}^{1+2n}$ endowed with the traslation operation $\circ$ defined in \eqref{eq:group-law}
      %and the family of dilations \eqref{eq:dil2}, 
      and we recover the estimate
   \begin{eqnarray} \label{eq:Cac1p}
			&&  \int_{B_{\varrho^{1+sp}}\times B_\varrho} \omega^2(\tau_2,x,v)  \dx\dv  
			+ \int_{-r^{sp}}^{\tau_2}\int _{B_{r^{1+sp}}} \mathcal{E}(f, \omega \varphi^p)\dt\dx\notag\\*
		  && \qquad\quad \leq  \frac{c\langle v_{\rm o}\rangle}{(r-\varrho)^{1+sp}} \int_{Q_r }\omega^2(t,x,v) \, \dt\dx\dv.
\end{eqnarray} 
which is the analogous to \eqref{eq:Cac1}.
    Here, we only treat the estimates regarding the third term on the left-hand side involving the energy of 
    the equation $\mathcal{E}_p(f, \omega \varphi^p )$, where $\omega = (f - \kappa)_+$ as in the proof of Lemma \ref{thm ccp}, and splitting it we obtain
    \begin{align*}
    & \int _{Q_r}
    	\int _{B_r} \snr{f(t,x,v)-f(t,x,w)}^{p-2}(f(t,x,w)-f(t,x,v)) \times \\
	& \hspace{1,5cm} \times ( \omega(t,x,v)\varphi^p(t,x,v) 
    	- \omega(t,x,w)\varphi^p(x,w) )K(t,x,v,w)\dw \dt \dx \dv \nonumber
    	\\*[0.5ex]
    	& + \, 2\int _{Q_r} \int_{\er^n \setminus B_r} 
       \snr{f(t,x,v)-f(t,x,w)}^{p-2}(f(t,x,w)-f(t,x,v)) \times \\
       & \hspace{5,2cm} \times \omega(t,x,v) \varphi^p(x,v) K(t,x,v,w) \dw \dt \dx \dv
    	\nonumber\\*[0.5ex]
    	&=: J_{2,1} + J_{2,2},
    \end{align*}
    We begin by estimating the term~$J_{2,1}$. 
   Firstly, if~$f (t,x,v) \geq f (t,x,w)$, then {
    \begin{align*}
    	&\snr{f (t,x,v) -f (t,x,w)}^{p-2}
	(f (t,x,v) -f (t,x,w))\\
	&\hspace{4cm} (\omega(t,x,v)\F^p(x,v) 
    	- \omega(t,x,w)\F^p (x,w) ) \\*[0.5ex]
    	& \,  = \big(f (t,x,v) -f(t,x,w)\big)^{p-1}\big(\omega(t,x,v)\F^p (x,v)  - \omega(t,x,v)\F^p (x,w) \big)\\*[0.5ex]
    	&  \, \geq
    	\begin{cases}
    		\big| \omega(t,x,v) -\omega(t,x,w)\big|^{p-1}
		 \times \\ 
		\hspace{0,6cm} \times \big( \omega(t,x,v)\F^p (x,v)  - \omega(t,x,w)\F^p (x,w) \big) 
		  \, \,  & \text{if} \, \, \, \ f (t,x,v) ,\, f (t,x,w)  >\kappa,\\
    		\big|\omega(t,x,v)\big|^{p-1} \omega(t,x,v) \varphi^p(x,v) \quad 
			\hspace{0,4cm} &\text{if} \, \,\, \ f (t,x,v) >\kappa \geq f (t,x,w),\\
    		0 \quad \hspace{5,2cm} & \text{otherwise,}  
    	\end{cases}\\*[0.5ex]
    	& \, \geq\ \big(\omega(t,x,v) -\omega(t,x,w)\big)^{p-1}
		\big( \omega(t,x,v)\F^p (x,v)  - \omega(t,x,w)\F^p (x,w) \big)\,,
    \end{align*} }
   which yields
     \begin{align*}
    & \snr{f (t,x,v) -f (t,x,w)}^{p-2}(f (t,x,v) -f (t,x,w)) \\
     		&\qquad \qquad (\omega(t,x,v) \F^p(x,v) - \omega(t,x,w) \F^p (x,w) ) K(t,x,v,w)\\
     \qquad  &\geq \big | \omega(t,x,v) -\omega(t,x,w)\big|^{p-2}
      	\big( \omega(t,x,v) -\omega(t,x,w)\big) \\
	&\qquad \qquad  \big( \omega(t,x,v) \F^p (x,v)  - \omega(t,x,w)\F^p (x,w) \big)K(t,x,v,w).
     \end{align*}
    If the opposite holds true, i.~\!e.~$f (t,x,v) \leq f (t,x,w)$,
    then we exchange the roles of~$v$ and~$w$ and repeat the computations above.
   
    Furthermore, under the assumptions~$\omega(t,x,v) \geq \omega(t,x,w)$ 
    and~$\F(x,w) \geq\F(x,v)$, by applying \cite[Lemma~3.1]{DKP14} we obtain
    \[
	    	(1-c_p\eps)\F^p(x,w) -(1 + c_p\eps)\eps^{1-p}|\F(x,v)-\F(x,w)|^p \leq \F^p(x,v). 
    \]
    Then, by choosing 
    \[
   	 \eps := \frac{1}{\max(1,2c_p)}
	 \frac{\omega(t,x,v)-\omega(t,x,w)}{\omega(t,x,v)} \in (0,1],
    \]
    we get
    \begin{eqnarray*} 
    &&\big(\omega(t,x,v)-\omega(t,x,w)\big)^{p-1}\omega(t,x,v)\F^p(x,v) \\
    && \geq  \big(\omega(t,x,v)-\omega(t,x,w)\big)^{p-1}
    		\omega(t,x,v) \max\{\F(x,v),\F(x,w)\}^p\\*[0.5ex]
    && \qquad -\frac{1}{2}\big(\omega(t,x,v)-\omega(t,x,w)\big)^p
    		\max\{\F(x,v),\F(x,w)\}^p\nonumber\\*[0.5ex]
    &&\qquad -c\max\{\omega(t,x,v),\omega(t,x,w)\}^p\snr{\F(x,v)-\F(x,w)}^p. 
    \end{eqnarray*}
 Then, one easily observes that the estimates above trivially hold 
 when~$0=\omega(t,x,v) = \omega(t,x,w)$, 
 or~$\omega(t,x,v) \geq\omega(t,x,w)$ and~$\F(x,v) \geq\F(x,w)$. 
 
 Hence, by only keeping the assumption~$\omega(t,x,v) \geq \omega(t,x,w)$, 
 we get
  \begin{eqnarray*} 
    &&\big(\omega(t,x,v)-\omega(t,x,w)\big)^{p-1}
    	\Big( \omega(t,x,v)\F^p(x,v) -  \omega(t,x,w)\F^p(x,w)  \Big )\\
    && \geq \big(\omega(t,x,v)-\omega(t,x,w)\big)^{p} \max\{\F(x,v),\F(x,w)\}^p\\*[0.5ex]
    &&\qquad -\frac{1}{2}\big(\omega(t,x,v)-\omega(t,x,w)\big)^p
    		\max\{\F(x,v),\F(x,w)\}^p\nonumber\\*[0.5ex]
    &&\qquad -c\max\{\omega(t,x,v),\omega(t,x,w)\}^p\snr{\F(x,v)-\F(x,w)}^p  \\
    &&\geq  \frac{1}{2}\big(\omega(t,x,v)-\omega(t,x,w)\big)^p
    		\max\{\F(x,v),\F(x,w)\}^p\\
	&&\qquad -c\max\{\omega(t,x,v),\omega(t,x,w)\}^p\snr{\F(x,v)-\F(x,w)}^p 	.
    \end{eqnarray*}
 Now, recalling that $\omega(t,x,v) \geq \omega(t,x,w)$, we rewrite the above  
estimate as follows
 \begin{eqnarray*} 
    &&\big|\omega(t,x,v)-\omega(t,x,w)\big|^{p-1}
    	\Big( \omega(t,x,v)\F^p(x,v) -  \omega(t,x,w)\F^p(x,w)  \Big )\\
    &&\geq  \frac{1}{2}\big|\omega(t,x,v)-\omega(t,x,w)\big|^p
    		\max\{\F(x,v),\F(x,w)\}^p\\
	&&\qquad -c\max\{\omega(t,x,v),\omega(t,x,w)\}^p\snr{\F(x,v)-\F(x,w)}^p ,
    \end{eqnarray*}
 	which is symmetric with respect to $v,w$. Hence, the above inequality holds true for every $v,w \in B_r$
	by exchanging the roles of~$v$ and~$w$.
  Finally, observing that
  \begin{eqnarray*}
 	&&\snr{\omega(t,x,v) \F(x,v)-\omega(t,x,w) \F(x,w)}^p \\ 
	&& \hspace{2cm} \leq  2^{p-1}
  	\snr{\omega(t,x,v)-\omega(t,x,w)} \max\{ \F(x,v),\F(x,w) \}^p\\*[0.5ex]
   && \hspace{2,5cm} + \, 2^{p-1}\snr{\F(x,v)-\F(x,w)}^p \max \{ \omega(t,x,v),\omega(t,x,w) \}^p,
  \end{eqnarray*}
   we conclude 
     \begin{align} \label{eq:J11-final_2}
     		J_{2,1}  &\geq  c  \int_{U_r} [\omega\varphi]^p_{W^{s,p}(B_r)} \dt \dx \nonumber \\
    		& \quad - c \int_{U_r}  \iint_{B_{r}\times B_{r}}\tfrac{\max\{ \omega(t,x,v),\omega(t,x,w) \}^p\snr{\F(x,v)-\F(x,w)}^p}{\snr{v-w}^{n+sp}} \dv \dw \dt \dx \nonumber \\
		&\geq c \int_{U_r} [\omega\varphi]^p_{W^{s,p}(B_r)} \dt \dx  - \frac{c\,r^{p(1-s)}}{(r-\rr)^p}\int_{Q_r}\omega^p(t,x,v)\dt\dx\dv,
    \end{align}
     where, once again we reasoned as in \eqref{eq:max} by symmetry of the Gagliardo kernel. 
     Now, we deal with the nonlocal term in~$J_{2,2}$. Note that
    \begin{eqnarray*}
    &&	\snr{f (t,x,v)-f (t,x,w)}^{p-2}(f (t,x,v)-f (t,x,w))\omega(t,x,v)  \\*[0.5ex]
	&&\quad \geq -\big(f (t,x,w)-f (t,x,v)\big)_+^{p-1}\big(f (t,x,v)-\kappa\big)_+  \\*[0.5ex]
    &&\quad \geq -\big(f (t,x,w)-\kappa\big)_+^{p-1}\big(f (t,x,v)-\kappa\big)_+  = -\omega^{p-1}(t,x,w)\omega(t,x,v).
    \end{eqnarray*}
   From this, we conclude that~$J_{2,2}$ can be treated as follows:
   \begin{align}\label{eq:J21_2}
   	J_{2,2} 	&\geq -c  \int_{Q_r} \int_{\er^n \setminus B_r}
   	\frac{\omega^{p-1}(t,x,w)\omega(t,x,v)\F^p(x,v)}{\snr{v-w}^{n+sp}} \dw\dt\dx\dv \\*[0.5ex]
   	&\geq -c\int_{Q_r} \omega(t,x,v)\F^p(x,v)\Biggl(\sup_{v \in \,\text{supp}\,\F}\int_{\er^n \setminus B_{r}}\frac{\omega^{p-1}(t,x,w)\dw}{\snr{v-w}^{n+sp}}\Biggr) \dt\dx\dv. 
	\nonumber
   \end{align}
  Hence, combining~\eqref{eq:J11-final_2} and~\eqref{eq:J21_2}, it yields that
  \begin{align}\label{eq:I22}
      J_{2} & \geq  c  \int_{U_r}[\omega\varphi]^p_{W^{s,p}(B_r)} \dt \dx- \frac{c\,r^{p(1-s)}}{(r-\rr)^p}\int_{Q_r}\omega^p(t,x,v)\dv\dx\dt\\*[0.5ex]
          &  -c  \int _{Q_r}  \omega(t,x,v)\F^p(x,v)\Biggl(\sup_{v \in \,\text{supp}\,\F}\int_{\er^n \setminus B_{r}}\frac{\omega^{p-1}(t,x,w)\dw}{\snr{v-w}^{n+sp}}\Biggr) \dt\dx\dv.
          \nonumber
  \end{align}
  All in all, by combining~\eqref{eq:Cac1p} and~\eqref{eq:I22} we obtain the desired result up to proceeding with the same estimates as at the end of Lemma \ref{thm ccp},
   for some constant $c$ depending only on $p$ and on the kernel constant $\Lambda$.
\end{proof}		

 Nevertheless, one has to face great difficulties when studying the weak regularity theory for \eqref{pbm4}, being the lack of a proper gain of integrability one  of the greater ones. 
 Indeed, the study of kinetic Sobolev spaces of the type $\W$ has not yet  been fully addressed, 
 despite some very recent work 
 \cite{PasPe}. Hence, there is no available proof  of any Sobolev-type embeddings in this specific framework, nor any information on the existence a fundamental 
 solution for the operator when $p\neq 2$. All in all,
 it appears it is not easy to overcome the obstacles posed by both the lack of ellipticity and by the $p$-growth of the integral diffusion by replacing the result of
Bouchut \cite{Bou02}, or the use of the fundamental solution constructed by Kolmogorov. 
Indeed, not even the apriori boundedness of weak solutions, such as in \cite{Sto19}, does lead
 to prove an explicit interpolative result in the fashion of Theorem \ref{loc bdd} without assuming, for instance, a suitable adaptation of the boundedness away from the vacuum assumption.

\end{document}